\documentstyle[12pt,amsfonts]{article}
\title{Expansiveness of algebraic actions on connected groups}
\author{Siddhartha Bhattacharya}
\begin{document}
\date{}
\newcommand{\Z}{{\mathbb{Z}}}
\newcommand{\R}{{\mathbb{R}}}
\newcommand{\Q}{{\mathbb{Q}}}
\newcommand{\C}{{\mathbb{C}}}
\newcommand{\G}{\Gamma}
\newcommand{\g}{\gamma} 
\newtheorem{theorem}{Theorem}[section]
\newtheorem{proposition}[theorem]{Proposition}
\newtheorem{lemma}[theorem]{Lemma}
\newtheorem{definition}[theorem]{Definition}
\newtheorem{corollary}[theorem]{Corollary}
\newtheorem{construction}[theorem]{Construction}
\def\rem{\refstepcounter{theorem}\paragraph{Remark \thetheorem}}
\def\rems{\refstepcounter{theorem}\paragraph{Remarks \thetheorem}}
\def\proof{\paragraph{Proof}}
\def\example{\refstepcounter{theorem}\paragraph{Example \thetheorem}}
\def\defin{\paragraph{Definition}}
\def\inter#1{\subsection*{\hbox{}\hfill{\normalsize\sl #1}\hfill\hbox{}}}
\maketitle
\section{Introduction}
Let $(X,d)$ be a metric space and $\rho$ be a continuous action of a
semigroup
$\Gamma$ on $X$. Then $(X,\rho)$ is said to be ${\it expansive \/}$
 if there exists an $\epsilon > 0$ such that for any two distinct points
$x,y$ in $X$,
$$ \sup_{\gamma\in\Gamma}\ d(\rho(\gamma)(x), \rho(\gamma)(y))
 \ge \epsilon .$$
Any such $\epsilon$ is called an expansive constant of $(X,\rho)$.
It is easy to check that when $X$ is compact,
the notion of expansiveness is independent
of the metric $d$.

If $X$ is a metrizable topological group and $\rho$ is an  action
of a semigroup $\Gamma$ on $X$ by continuous endomorphisms of $X$, 
then $\rho$ is said to be
{\it expansive} if there exists a neighborhood $U$ of the identity
in $X$ such that 
$${\bigcap}_{\gamma} \rho(\gamma)^{-1}(U) = \{ e\} .$$
Any such neighborhood is called an {\it expansive neighborhood}
of the identity in $X$. An automorphism  $\tau$ of $X$ is
said to be expansive  if the cyclic
group  generated by $\tau$ acts expansively on
$X$.
It is easy to check that when $X$ is compact and $\rho$ is an
endomorphism action of a semigroup $\Gamma$, these two
notions of expansiveness coincide. 

The notion of expansiveness plays an important role in the study of 
dynamical systems in general and endomorphism actions in particular.
Earlier, expansiveness of automorphisms  of  connected groups was  
studied by
several authors (cf. \cite{A}, \cite{A-D}, \cite{E1}, \cite{E2} 
 and  \cite{Law}). 
A complete description of expansive automorphisms
on compact connected groups was obtained in \cite{Law}. 
In particular, it was shown 
that every  compact connected group which admits an 
expansive automorphism is abelian and finite-dimensional.
Also in \cite{La} it was proved that if a compact
 connected group $X$ admits an expansive endomorphism action
then it is abelian.

In recent years dynamics of endomorphism actions of discrete
semigroups 
on compact abelian groups has been extensively studied using
techniques from commutative algebra (cf. \cite{Sc2}). 
Among other things, a complete characterization of 
expansiveness has been obtained when $X$ is zero-dimensional 
( cf. \cite{Ki-Sc}, Theorem 5.2) or $\Gamma = {\mathbb{Z}}^{d}$ for
 some $d\ge 1$
 (cf. \cite{Sc1}, Theorem 3.9). For further study of expansiveness and 
it's relation with other dynamical properties the reader is referred
to \cite{B-L}, \cite{Ka-Sc} and \cite{L-S}.

\medskip
In this paper we consider expansive  endomorphism 
 actions of arbitrary  semigroups
on connected metrizable topological groups. For any such action
$(G,\rho)$ we give  necessary and sufficient conditions for 
expansiveness of $\rho$, provided $G$ is either a Lie group or a
 compact finite-dimensional group. 

\medskip
This paper is organized as follows. In section 2 we prove some
elementary results about expansiveness of endomorphism actions on
finite-dimensional vector spaces, considered as abelian topological
groups under addition.
In Section 3 we consider endomorphism
actions on connected Lie groups. If $G$ is a connected Lie group then by 
$L(G)$ we denote the Lie algebra of $G$. If $\rho$ is an
endomorphism action of a semigroup $\Gamma$ on $G$ then by 
$\rho_{e}$ we denote the 
induced $\Gamma$-action on $L(G)$ defined by
$$ \rho_{e}(\gamma) = {\rm d}\rho(\gamma)|_{e}\ \ \forall \gamma \in
\Gamma.$$
 We prove the following.

\medskip
\noindent
{\bf Theorem A.}
{\it  
Let $G$ be a connected Lie group, $\Gamma$ be a semigroup and
$\rho$ be an endomorphism action of $\Gamma$ on $G$. Then
 $(G,\rho)$ is
expansive if and only if for all non-zero $v$ in $L(G)$,             
 the $\rho_{e}$-orbit of $v$ is unbounded. \/}

\bigskip
 In  the special case when 
 $\Gamma$ is either abelian or a virtually nilpotent
 group, applying the above theorem we give a necessary and sufficient 
condition for  
expansiveness of $(G,\rho)$ in terms of the generalized weights of $\rho_{e}$.
 As a consequence we
 prove that  if $\Gamma$ is a virtually nilpotent group and 
$(G,\rho)$ is expansive then for some $\gamma_{0}$ in $\Gamma$, 
$\rho(\gamma_{0})$ is 
an expansive automorphism of $G$. 

\medskip
In section 4 we consider endomorphism actions on compact connected
finite-dimensional abelian groups. For any such group $G$, by
$\widehat{G}$ we denote the Pontryagin dual of $G$ and
by $L(G)$ we denote the vector space of all homomorphisms from
$ \widehat{G}$ to ${\mathbb R}$ under pointwise
addition and scalar multiplication. From the duality theory of compact
 abelian groups it follows that $L(G)$ is finite-dimensional. 
If $\rho$ is an
endomorphism action of a semigroup $\Gamma$ on $G$ then by 
$\rho_{e}$ we denote the 
induced $\Gamma$-action on $L(G)$ defined by
$$ \rho_{e}(\gamma)(p)(\chi) = p(\chi\circ \rho(\gamma))\ \ \ \forall
p\in L(G).$$
 Also, for such an
action $\rho$, by $\widehat{\rho}$ we denote the induced
$\Gamma$-action on $\widehat{G}$. We note that  
$\widehat{G}$ can be
realised as a module over ${\mathbb Z}(\Gamma)$, the group-ring of 
$\Gamma$, via the action $\widehat{\rho}$.

\medskip
\noindent
 We prove the following.

\medskip
\noindent
{\bf Theorem B.}
{\it 
Let $\Gamma$ be a semigroup and
$\rho$ be an endomorphism action of $\Gamma$ on
 a compact connected finite-dimensional abelian group $G$.
Then $(G,\rho)$ is expansive 
if and only if the following two conditions are satisfied.

\medskip
a) $\widehat{G}$ is finitely generated as a ${\mathbb Z}[\Gamma]$-module.

\medskip
b) For every non-zero $v$ in $L(G)$, the $\rho_{e}$-orbit of
$v$ is unbounded. \/}
\rem
Since non-abelian  compact connected groups do not admit  expansive
endomorphism actions  (cf. \cite{La}),
  Theorem B characterizes 
expansive endomorphism actions 
on arbitrary compact connected finite-dimensional groups.
 
\section{Preliminaries}
Throughout this paper $\Gamma$ will denote a discrete semigroup i.e. a
set
with an associative multiplication law. If $\Gamma_{1}, \Gamma_{2}$ are
semigroups then a map $f$ from $\Gamma_{1}$ to $\Gamma_{2}$ is said to be a
semigroup homomorphism if for all $\gamma_{1},\gamma_{2}$ in $\Gamma_{1}$,
$f(\gamma_{1}\cdot\gamma_{2}) = f(\gamma_{1})\cdot
f(\gamma_{2})$. If $G$ is a topological group then by
${\rm End}(G)$ we denote semigroup of all continuous endomorphisms of
$G$, with
composition as the product law.

Suppose $G$ is a topological group and $\rho : \Gamma\rightarrow {\rm
End}(G)$
is an endomorphism action of a discrete semigroup $\Gamma$ on $G$.
Then it is easy to see that $\Gamma_{1} = \rho(\Gamma)\cup\{{\rm
Id}\}$ is a subsemigroup of ${\rm End}(G)$ and $\rho$ is 
expansive if and only if the natural action of $\Gamma_{1}$ 
 on $G$ is expansive. 
Henceforth we will restrict the discussion of endomorphism actions to
the case when 
$\Gamma$ contains an identity element $e$ and $\rho(e)\in {\rm End}(G)$
is the identity automorphism.

\medskip
\noindent
 If ${\bf C}$ denotes the semigroup of complex numbers
under multiplication then any  
semigroup homomorphism $\lambda : \Gamma \rightarrow {\bf C}$ is said
 to be a ${\it character \/}$ of $\Gamma$.
If $\rho : \Gamma\rightarrow {\rm End}(V)$ is an endomorphism action
of $\Gamma$ on  a finite-dimensional vector space $V$ over $\C$,  
then $\rho$ is said to be {\it reducible \/} if there exists
characters  $\lambda_{1},\ldots ,\lambda_{k}$
and nontrivial $\rho$-invariant subspaces
$V_{1},\ldots ,V_{k}\subset V$
satisfying the following conditions.
\begin{enumerate}
\item{ $ V = V_{1}\oplus \cdots \oplus V_{k}$.}
\item{
  For  $i = 1,\ldots ,k$, there exists a basis $B$ of
$V_{i}$ such that 
$\rho(\gamma)|_{V_{i}} - \lambda_{i}(\gamma)I$ is strictly
upper triangular with respect to $B$, for all $\gamma$ in $\Gamma$.}
\end{enumerate}
The characters $\lambda_{1},\ldots ,\lambda_{k}$ are said to be the
{\it generalized weights \/} of $\rho$ and for all $i$, 
$V_{i}$ is said to be the
{\it generalized weight space \/} corresponding to $\lambda_{i}$.
\example
If  $\rho(\Gamma)$ is an  abelian subsemigroup
 ${\rm End}(V)$ then $\rho$ is reducible (see \cite{Ja}, pp 134).
 This is also true when $\rho(\Gamma)$  is
contained in a connected nilpotent subgroup $H$ of $GL(V)$. This
can be  
seen by looking at the Lie algebra homomorphism
$i : L(H) \rightarrow {\mbox{End}}(V)$ and
using certain fundamental facts about representations of nilpotent
Lie algebras (see \cite{Kn}, Proposition 2.4 and Theorem 1.35).

\medskip
The following proposition gives a necessary and sufficient condition
for expansiveness of endomorphism actions on a finite-dimensional
vector space, considered as an abelian group under addition.
\begin{proposition}\label{p1} 
Let $\Gamma$ be a semigroup and $\rho : \Gamma\rightarrow {\rm
End}(V)$ be an endomorphism action of $\Gamma$ on a 
 finite-dimensional vector space $V$ over
${\mathbb R}$ or ${\mathbb C}$. Then we have the following :

\medskip
\noindent
a) $\rho$ is expansive if and only if 
for  all non-zero $v$ in $V$,  the 
$\Gamma$-orbit
of $v$ is unbounded. \\
b) If $V$ is a real vector space then
$\rho$ is expansive if and only if the induced endomorphism action of
$\Gamma$ on $V\otimes {\mathbb C}$ is expansive.
\end{proposition}
\proof
For any $v$ in $V$ let $O_{v}\subset V$ denote the $\rho$-orbit of $v$. 
Suppose that there is a non-zero vector $v$ in $V$ with bounded
$\rho$-orbit.  
Let $U$ be any neighborhood of $0$ in $V$.
Then for $\epsilon > 0$ sufficiently small, 
$O_{\epsilon v} = \epsilon O_{v}$ is contained in $U$. 
This implies that $\rho$
is not expansive. On the other hand if $O_{v}$ is unbounded for all non-zero
$v\in V$,
then it is easy to see that any bounded neighborhood of $0$ is 
an expansive  neighborhood. This proves a). Part b) is an
immediate consequence of a). 
$\hfill \Box$ 
\begin{proposition}\label{p2}
Let $\Gamma$ be a semigroup, $V$ be a   
 finite-dimensional vector space over ${\mathbb C}$ and
 $\rho : \Gamma\rightarrow {\rm End}(V)$ be a reducible endomorphism
action with  generalized weights $\lambda_{1},\ldots ,\lambda_{k}$.
Then $\rho$ is expansive  if and
only if for each $i$,
 the image of $\lambda_{i}$ is an unbounded subset of ${\mathbb C}$.
\end{proposition}
\proof
For $i = 1,\ldots ,k$ let $V_{i}$ be the generalized weight space corresponding to
$\lambda_{i}$. Since each $V_{i}$ is $\rho$-invariant,
 from the previous proposition  it follows that $\rho$ is
expansive
if and only if $\rho|_{V_{i}}$ is expansive for all $i$.
Therefore without loss of generality we may assume that
$k = 1$. In that case there exists a character $\lambda$ of $\Gamma$
such that after suitable identifications we have :
\\
a) $V = {\mathbb C}^{n}$ for some $n > 0$ and
\\
b) For all $\gamma\in\Gamma$, $\rho(\gamma)$ can be written as,
$$ \rho(\gamma)\  =\   
\left(  
\begin{array}{ccccc} 
 \lambda(\gamma) &  &  &  &  \\
  & \ddots & & * & \\
 & & \lambda(\gamma)& & \\
 & 0 & & \ddots  & \\
 & & & & \lambda(\gamma)   
\end{array}
\right) . $$
Suppose there exists a non-zero vector $v = (v_{1}\ldots ,v_{n})^{t}$
in $\C^{n}$ such that the $\rho$-orbit of $v$ is bounded. Let 
$j \le n$ be the largest positive integer such that $v_{j} \ne 0$.
Then it is easy to see that
$$ \{ \rho(\gamma)(v)_{j}\ | \ \gamma\in \Gamma\} = 
\{ \lambda(\gamma)v_{j}\ | \ \gamma\in \Gamma\}. $$
Since  the $\rho$-orbit of $v$ is bounded, this implies that
the image of $\lambda$ is bounded subset of $\C$.
Conversely, it is easy to see that if $v\in\C^{n}$ denotes
the vector $(1,0,\ldots ,0)^{t}$ then
$$ \{ \rho(\gamma)(v)\ | \ \gamma\in \Gamma\} = 
 \{ (\lambda(\gamma),0,\ldots ,0)^{t}\ | \ \gamma\in \Gamma\}.$$
Therefore if the image of $\lambda$ is bounded then there exists
a non-zero vector in $V$ with bounded $\rho$-orbit. 
$\hfill \Box$

\medskip
In the special case when $\Gamma$ is isomorphic to 
$ \Z^{+}$ or $\Z$ we obtain the following consequence of the 
above proposition. This was first proved in  \cite{E2}.
\label{red}
\begin{proposition}\label{p3}
Let $V$ be a  finite-dimensional vector space over $\mathbb{R}$ or
$\mathbb{C}$ and let
$T$ be an element of End$(V)$. Then the cyclic semigroup generated by
$T$ acts expansively
on $V$ if and
only if the spectrum of $T$ does not intersect the closed unit disk.
If $T$ is invertible then the cyclic subgroup generated by $T$
acts expansively
on $V$ if and
only if the spectrum of $T$ does not intersect ${{\mathbb{S}}^{1}}$.
\end{proposition}
\proof 
Replacing $V$ by $V\otimes {\mathbb{C}}$ if necessary, we may assume
that $V$ is a complex vector space.
Let $\Gamma = \{T^{n}\ |\ n\ge 0\}$ denote the
cyclic semigroup generated by $T$ and let $\rho$ denote the
natural action of $\Gamma$ on $V$.   
Since $\Gamma$ is abelian, $\rho$ is reducible. Furthermore if 
$\lambda_{1},\ldots ,\lambda_{k}$ denote the generalized 
weights of $\rho$ then
$\lambda_{1}(T),\ldots ,\lambda_{k}(T)$
are the eigenvalues of $T$. Applying the previous proposition
we conclude that $\rho$ is expansive if and only if
the spectrum of $T$ does not intersect the closed unit disk.
 The second assertion follows from a similar argument.
$\hfill \Box$
\section{Expansive actions on Lie groups}
In this section we prove Theorem A and it's consequences. Throughout the
section $G$ will denote a connected Lie group and $L(G)$ will denote the
Lie algebra of $G$. The standard exponential map from $L(G)$ to $G$ will
be denoted by exp.

\medskip
We
begin with the
following lemma.
\begin{lemma}\label{key}
Let  $V$ be a finite-dimensional vector space over
${\mathbb R}$ and let $\rho$ be an
expansive endomorphism action of a semigroup $\Gamma$ on $V$.
Then there exists a finitely generated subsemigroup $\Gamma_{0} 
\subset \Gamma$ 
which acts expansively on $V$.
\end{lemma}
\proof
Suppose this is not the case.
We fix a norm $||.||$ on $V$. Let $S\subset V$ denote the set defined
by
 $ S = \{ v\in{V}\ |\ ||v|| = 1\}$.
For every
finitely generated subsemigroup $\Lambda \subset \Gamma$ we define
a closed subset $S({\Lambda})\subset S$ by
$$ S(\Lambda) = \{ v\in{S}\ |\ ||\rho(\gamma)(v)|| \le 1 \ \forall \gamma
\in{\Lambda}\}.$$
We claim that for any 
finitely generated subsemigroup $\Lambda \subset \Gamma$,
$ S(\Lambda)$ is non-empty.
To see this we define 
a $\Lambda$-invariant subspace $V(\Lambda)\subset V$ by
$$ V(\Lambda) = \{ v\in{V}\ |\ \Lambda-{\rm orbit\  of}\ 
 v\ {\rm is\  bounded}\}.$$ 
By Proposition \ref{p1}, $V(\Lambda)$ is a non-trivial subspace of $V$. 
Let $\Lambda_{b}\subset {\rm End}(V(\Lambda))$ be the image of
$\Lambda$ under the map $\gamma \mapsto \rho(\gamma)|_{V(\Lambda)}$.
Then $\overline{\Lambda_{b}}$ is a compact subsemigroup
of ${\rm End}(V(\Lambda))$.
We choose a non-zero $w\in{V(\Lambda)}$ and
define a continuous function 
$h : {\overline{\Lambda_{b}} }\rightarrow \R$ by
$h(\alpha) = || \alpha(w)||$ for all $\alpha$ in 
$\overline{\Lambda_{b}}$.
 Since $\overline{\Lambda_{b}}$
 is compact, there
exists a $\alpha_{0}\in{\overline{\Lambda_{b}}}$ such that
$h(\alpha_{0}) \ge h(\alpha)$ for all $\alpha$ 
in $\overline{\Lambda_{b}}$. 
 Putting
$v = \alpha_{0}(w)$ we see that
 for any $\alpha$ in $\overline{\Lambda_{b}}$,
$$ || \alpha(v) || = h(\alpha \alpha_{0}) \le h(\alpha_{0}) = ||v||.$$
 Since 
 $\overline{\Lambda_{b}}$ contains  Id,
$||\alpha_{0}(w)|| \ge ||w|| > 0$. Now it is easy to see that
 the unit vector $v/||v|| \in S(\Lambda)$, which proves the claim.

If
 $\Lambda_{1}, \Lambda_{2},\ldots ,\Lambda_{k}$ are
 finitely generated subsemigroups of $\Gamma$ and  $\Lambda$ is
a finitely generated subsemigroup of $\Gamma$ containing
$\Lambda_{1}, \Lambda_{2},\ldots ,\Lambda_{k}$, then
from the above claim it follows that
 $$ \emptyset \ne S({\Lambda})\subset S(\Lambda_{1})\cap\cdots\cap
S(\Lambda_{k}).$$
This shows that the collection 
$ \{ S( \Lambda )\}$ has the finite-intersection
property.
 Now from the compactness of $S$ it follows that 
$$ \bigcap_{\Lambda} S(\Lambda) \ne \emptyset .$$
Clearly for any vector $w$ which lies in the intersection, the
 $\Gamma$-orbit of $w$ is bounded, which contradicts our hypothesis.
$\hfill \Box$

\bigskip
{\bf Proof of Theorem A. }
Suppose $(G,\rho)$ is expansive. Let $U$ be an expansive neighborhood 
of the identity in $G$. We choose a neighborhood $V$ of $0$ in $L(G)$
 such that ${\rm exp}|_{V}$ is a homeomorphism and exp($V$) is contained in
$U$. Let $v$ be any non-zero element of $L(G)$ such that the 
$\rho_{e}$-orbit of $v$ is contained in $V$.
Since exp  is a
$\G$-equivariant map from $(L(G),\rho_{e})$ to $(G,\rho)$, 
it follows that the $\rho$-orbit of exp($v$) is contained in $U$.
Since $U$ is an expansive neighborhood 
of the identity  and ${\rm exp}|_{V}$ is a homeomorphism,
this implies that $v = 0$.  Hence $\rho_{e}$ is expansive. Now applying
Proposition \ref{p1} we see that for every 
non-zero $v$ in $L(G)$, the $\rho_{e}$-orbit of $v$
is unbounded.

Now  we will prove the converse. Applying Lemma \ref{key} and 
Proposition \ref{p1} we see 
that it is enough to consider the case when $\G$ is finitely
generated.  Let
$A$ be a finite set which generates $\G$.
We choose a neighborhood 
$V$ of $0$ in $L(G)$
such that ${\rm exp}|_{V}$ is a homeomorphism and $\overline{V}$ is
compact. We choose another neighborhood $U$ of $0$ such that
$$U\subset V\ \ {\rm and}\ \ \rho_{e}(\gamma)(U)\subset V
\ \forall \gamma \in A.$$ 
We claim that exp($U$) is an expansive neighborhood of $(G,\rho)$.
 Suppose this is not the case. Let $g\ne e$ be any element of $G$ such
that the $\rho$-orbit of $g$ is contained in ${\rm exp}(U)$. 
  We choose $v$ 
in $U$ such that ${\rm exp}(v) = g$. Let $O_{v}$ be
the $\rho_{e}$-orbit of $v$.  Since $v$ is non-zero,  $O_{v}$ is
unbounded.
Since $A$ generates $\G$ and
$v\in{U}$,
it follows that  there exists
$w$ in $O_{v}$  and $\gamma_{0}\in A$ such that
$$w\in U \mbox{ and } \rho_{e}(\gamma_{0})(w)\in V - U.$$ 
Since exp maps the $\rho_{e}$-orbit of $v$ onto the 
 $\rho$-orbit of exp($v$), this implies
that the $\rho$-orbit of ${\rm exp}(v)$ intersects  ${\rm exp}(V - U)$. 
Since ${\rm exp}(V - U)$ and ${\rm exp}(U)$ are disjoint, this
gives a contradiction. 
$\hfill \Box$

\medskip
We note the following consequence of
 Theorem A and Lemma \ref{key}.
\label{c1} 
\begin{corollary}\label{im}
Let $\G$ be a discrete semigroup, $G$ be a connected Lie group and $\rho$
be an expansive endomorphism action of $\G$ on $G$. 
Then there exists a finitely
generated subsemigroup $\G_{0}\subset \G$ such that the action of $\rho$,
restricted to $\G_{0}$, is expansive.
\end{corollary}  
The following corollary provides a rich source of expansive
endomorphism actions of non-abelian semigroups on tori.

\label{c2}
\begin{corollary}\label{tor}
Let $\Gamma$ be an infinite subsemigroup of $M(n,{\mathbb Z})$ which
acts irreducibly  on ${\mathbb R}^{n}$.
Then the induced $\Gamma$-action on 
${\mathbb T}^{n}\cong {\mathbb R}^{n}/ {\mathbb Z}^{n}$ is expansive.
\end{corollary}
\proof
Let $\rho$ denote the natural action of $\Gamma$ on ${\mathbb R}^{n}$
and $W$ denote the subspace of ${\mathbb R}^{n}$ consisting of all
points with bounded $\rho$-orbit. Since $\Gamma$ is an
infinite subsemigroup of $M(n,{\mathbb Z})$,
 it is a noncompact subset of 
of $M(n,{\mathbb R})$. Hence $W$ is a proper subspace of ${\mathbb R}^{n}$.
From irreducibility of $\rho$ we conclude that $W = \{0\}$. Now
applying Theorem A we see that the action of $\Gamma$ on ${\mathbb T}^{n}$
is expansive.
$\hfill \Box$
 \label{cy2}
\begin{corollary}\label{nil}
Let $G$ be a  connected Lie group and $\rho$ be an expansive automorphism
action of a  virtually nilpotent group $\Gamma$ on $G$.
 Then $\rho(\G )$ contains an expansive automorphism of $G$.
\end{corollary}
\proof
Since $\Gamma$ acts expansively on $G$, so does every finite-index
subgroup of $\Gamma$.
Therefore without loss of generality we may assume that
 $\Gamma$ is nilpotent.
First we will consider the case when $G$ is isomorphic to a 
finite-dimensional vector $V$ space over ${\mathbb C}$.
 Let $H\subset GL(V)$ be the Zariski-closure
of $\rho(\G )$ and $H_{0}$ be the connected component of $H$
which contains the identity. Since $H$ has finitely many components,
it follows that 
$$ \G_{0} = \{ \gamma\in \G\ |\ \rho(\gamma)\in H_{0}\}$$
is a finite-index subgroup of $\G$. Let $\rho_{0}$ denote the
$\Gamma_{0}$-action on $V$, induced by $\rho$. Since $\rho(\Gamma_{0})$
 is contained in a connected
nilpotent subgroup of $GL(V)$, $\rho_{0}$ 
is reducible. Since $\rho$ is expansive, so is $\rho_{0}$. Let
$\lambda_{1},\ldots ,\lambda_{k}$ be  the generalized weights  
  of $\rho_{0}$. 
 For each $\gamma$ in $\G_{0}$ we define 
$A_{\gamma}\subset \{1,\ldots ,k\}$ by
$$ A_{\gamma} = \{ j\ |\ |\lambda_{i}(\gamma)| \ne 1\}.$$
For any $\g\in\G_{0}$ let $n(\g)$ denote the cardinality of
$A_{\g}$.
We choose $\g_{0}\in \G_{0}$ such that the $n(\g_{0}) \ge n(\g)$ for all
$\g$ in $\G_{0}$. We claim that $n(\g_{0}) = k$ i.e.
$ |\lambda_{i}(\g_{0})| \ne 1$ for all $j = 1,\ldots ,k$.
Suppose this is not the case. We choose $i$ such that 
$ |\lambda_{i}(\g_{0})| = 1$.                                 
By Proposition \ref{p2} there exists $\g_{1}$ such that  
$ |\lambda_{i}(\g_{1})| \ne 1$. 
It is easy to see that for sufficiently large $m > 0$,
$$ |\lambda_{j}(\g_{0}^{m}\g_{1})| \ne 1 \ \forall j\in A_{\g_{0}}
\cup \{i\}.$$
 Since $n(\g_{0}) \ge n(\g)$ for all
$\g$ in $\G_{0}$, this gives a contradiction. 
Since the numbers $\lambda_{1}(\gamma_{0}),\ldots ,\lambda_{k}(\gamma_{0})$ are
the eigenvalues of $\rho(\gamma_{0})$, from the above
claim and
 Proposition \ref{p3}  we 
see that $\rho(\gamma_{0})$ is an expansive automorphism of $V$.

Now we will consider the general case. Let $\sigma$ be the 
$\Gamma$-action  on the complexified
Lie algebra $  L(G)\otimes {\mathbb C}$, induced by $\rho_{e}$.
 Applying Theorem A and
Proposition \ref{p1} we see that $\sigma$ is expansive.
By the previous argument there exists a $\gamma_{0}$ in $\Gamma$ such
that
$\sigma(\gamma_{0})$ is an expansive automorphism of $  L(G)\otimes
{\mathbb C}$.  Applying Theorem A and
Proposition \ref{p1} to the cyclic group generated by $\gamma_{0}$,
we conclude that $\rho(\gamma_{0})$ is an expansive automorphism of $G$.
$\hfill \Box$

\medskip
It is known that if a connected locally compact topological group
 admits an
expansive automorphism then it is nilpotent (see \cite{A}, \cite{Mu}). The
following corollary generalises this result for actions of virtually
nilpotent groups on connected Lie groups.
\begin{corollary}\label{nil2}
Let $G$ be a  connected Lie group which admits an expansive 
automorphism action of a virtually nilpotent group. 
Then $G$ is nilpotent.
\end{corollary}
\proof
Let $\rho$ be an automorphism action of a virtually nilpotent group
$\G$ on $G$.
 Then by the previous corollary there exists a $\gamma
\in\G$ such that $\rho(\gamma)$ is an expansive automorphism of $G$.
 Applying Theorem A and Proposition \ref{p3} we see that the 
spectrum of $d\gamma \in {\rm Aut}(L(G))$ does not intersect the unit
circle. 
It is known that  if a finite-dimensional Lie algebra over ${\mathbb R}$
 admits 
a hyperbolic  automorphism then it is nilpotent. Hence $G$ is nilpotent.
$\hfill \Box$

\medskip
We conclude this section with an example showing that Corollary \ref{nil}
 does not hold for actions of arbitrary discrete groups.
\example
Fix $n\ge 3$ and define a  subgroup $\Gamma$ of $GL(n,{\mathbb R})$
 by
$$ \Gamma =  \left\{
\left(  
\begin{array}{cc} 
 A & b  \\
 0 & 1  \\   
\end{array}
\right)
\ \mid\  A\in{GL(n-1,{\mathbb Z})}, b\in{{{\mathbb Z}}^{n-1}}\
 \right\}  $$
 It is easy to see that
for any $x = (x_{1},\ldots ,x_{n})$
in $ {\mathbb R}^{n}$, the
$\Gamma$-orbit of $x$ is given by the set
$$\{ Ay + x_{n}b\ |\ A\in{GL(n-1,{\mathbb Z})}, b\in{{{\mathbb Z}}^{n-1}}\
 \}\, ;\ y = (x_{1},\ldots ,x_{n-1}).$$
Hence for every non-zero $x$ in $ {\mathbb R}^{n}$, the
$\Gamma$-orbit of $x$ is
unbounded. Since  the spectrum of $\gamma$
contains 1 for every $\gamma$ in $\G$, applying 
Proposition \ref{p3} we see that $\Gamma$ does
not contain any expansive automorphism of $ {\mathbb R}^{n}$.
\section{Expansive actions on solenoids}
In this section we consider endomorphism actions on solenoids (compact
 connected finite-dimensional abelian groups).
We freely use various results from duality theory of locally compact
abelian groups; the reader is referred to \cite{P} for details.

For any locally compact
abelian group $G$, we denote
by $\widehat G$ the dual group of $G$.
Recall that for a
compact connected  abelian group $G$, we denote
by $L(G)$ the vector space consisting of all homomorphisms
from $\widehat G$ to ${\mathbb R}$, under pointwise addition and
scalar multiplication.
It is known that if $G$ is a solenoid, then $\widehat G$ is a torsion free
discrete abelian group of finite rank. Hence for any solenoid $G$,
$L(G)$ is 
finite-dimensional.
 
If  $\rho$ is an
endomorphism action of a semigroup $\Gamma$ on a compact connected
abelian group $G$ then by 
$\widehat{\rho}$ we denote the induced
$\Gamma$-action on $\widehat{G}$ and by 
$\rho_{e}$ we denote the 
induced $\Gamma$-action on $L(G)$ defined by
$$ \rho_{e}(\gamma)(p)(\chi) = p(\chi\circ \rho(\gamma)).$$
\example
Let $\G$ be a semigroup, ${\widehat{\rho}}:\G\rightarrow M(n,\Q)$ be 
a semigroup homomorphism for some $n\ge 1$  
and $H$ be a $\G$-invariant subgroup of $\Q^{n}$.
 Then the  action of $\G$ on $H$ induces an 
endomorphism action 
$\rho$ of $\G$ on $G = \widehat{H}$. If $H$ contains $\Z^{n}$ then
 $L(G)$ can be identified with $\R^{n}$ and $\rho_{e}$ can be
identified with the  adjoint of the  $\G$-action  on  $\R^{n}$ induced
by $\widehat{\rho}$.
From duality theory of compact abelian groups it follows that
any endomorphism action of a  discrete semigroup on a solenoid
 can be identified with an action
of this form.

\medskip
This section is organized as follows. Throughout this section
$G$ will denote a solenoid and $\rho$ will denote an
endomorphism action of a discrete semigroup $\Gamma$ on $G$. 
In 4.1 we prove that the conditions a) and b), as stated in Theorem B are
necessary for expansiveness of $\rho$. In
4.2 we  show that  expansiveness of $\rho$ and  existence of non-trivial
bounded $\rho_{e}$-orbit in $L(G)$ can be characterized in terms of
suitably chosen metrics on $G$ and $L(G)$ respectively. In 4.3 we prove
Theorem B when $\Gamma$ is finitely generated. We complete the proof of
Theorem B in 4.4.
\subsection{Necessary conditions for expansiveness}
For any 
compact connected  abelian group $G$, 
we define a map $E$ from $L(G)$ to $G$ by the condition
$$(\phi \circ E )(p) = e^{2\pi i p(\phi )} \ \ \ \forall p\in{L(G)} ,
 \phi \in{\widehat G}.$$
Since for a fixed $p$ in $L(G)$, the map 
$\phi \rightarrow e^{2\pi i  p(\phi )}$ is a continuous homomorphism from
$\widehat{G}$ to ${{\mathbb{S}}^{1}}$, by the duality theorem the map $E$ is 
well defined and unique. From the uniqueness it follows that 
$E$ is a homomorphism from $L(G)$ (considered as an abelian group
under addition) to $G$. It is easy to check that the kernel of $E$ can
be identified with the set of all homomorphisms from $\widehat{G}$ to
${\mathbb Z}$, which is a discrete subgroup of $L(G)$.
\rem 
 Using the duality theorem we can realize $L(G)$ with the 
set of all one-parameter subgroups of $G$ and $E$ with the map
$\alpha \mapsto \alpha (1)$. In particular when $G$ is a
 torus, $L(G)$   
can be identified with ${\mathbb R}^{n}$, the Lie
algebra of $G$, and $E$ can be identified with the standard exponential
map. However in general the map  $E$ is not
surjective. In fact from a result of Dixmier (see \cite{D}) it follows that
if $G$ is metrizable then $E$ is surjective 
if and only if $\widehat{G}$ is a free abelian
group.
\defin
Let $\Gamma $ be a  semigroup and $\rho $ be an  endomorphism
action of $\Gamma $ on  a compact abelian group $G$.  
Then a set $A\subset\widehat{G}$ is said to be a {\it $\rho$-basis \/}
if $A$ generates $\widehat{G}$ as an abelian group 
 and $A = \bigcup_{\gamma}{\widehat{\rho}}(\gamma)(F)$ for some
finite set $F \subset \widehat{G}$. 

\medskip
Clearly $\widehat{G}$ admits a $\rho$-basis if and only if it is
finitely generated as a ${\mathbb Z}(\Gamma)$-module.

\medskip
The next two propositions show that for expansive endomorphism actions
the conditions a) and b), as stated in Theorem B are satisfied. The
results are known, we include the proofs for the sake of completeness. 
\label{cn1}
\begin{proposition}\label{basis}
 Let $\Gamma $ be a  semigroup and $\rho $ be an expansive endomorphism
action of $\Gamma $ on  a compact abelian group $G$.  
Then $\widehat{G}$ has a $\rho$-basis.
\end{proposition}  
\proof
Let $U$ be an expansive neighborhood of of $e$ in $G$. We choose a
finite set $F\subset {\widehat{G}}$ and $\epsilon > 0 $ such that
$$  \bigcap_{\chi\in{F}} \{ g\in G\ |\ |\chi(g) - 1 |
 < \epsilon\}\subset U.$$
We define $A \subset {\widehat{G}}$ by 
$ A = \bigcup_{\gamma}{\widehat{\rho}}(\gamma)(F)$.
Let ${\widehat{H}}\subset{\widehat{G}}$ be the subgroup generated by
$A$ and let $G^{'}\subset G$ be the subgroup consisting of all
$g$ in $G$ such that $\chi(g) = e$ for all $\chi$ in ${\widehat{H}}$.
Then for every $g\in G^{'}$, the $\rho$-orbit of $g$
is contained in $U$. Since $U$ is an expansive neighborhood of of $e$,
we conclude that $G^{'} = \{ e\}$. Now from duality theory of compact
abelian groups it follows that 
${\widehat{H}} = {\widehat{G}}$ (see \cite{P}, Theorem 53). Hence $A$ is a
$\rho$-basis of $\widehat{G}$. 
$\hfill \Box$
\label{cn2}
\begin{proposition}\label{div}
(Also see \cite{B}, proof of Corollary 1)
Let $\Gamma $ be a  semigroup and $\rho $ be an expansive endomorphism
action of $\Gamma $ on  a solenoid $G$.  
 Then for every non-zero point
$p\in{L(G)}$, the orbit of $p$ under $\rho_{e}$ is unbounded.
\end{proposition}
\proof
Let $U$ be an expansive neighborhood
of $e$ in $G$.  Since $E$ is a $\G$-equivariant map
from $(L(G),\rho_{e})$ to $(G,\rho)$, it follows that
$$
 \bigcap_{\gamma\in\G}\rho_{e}(\gamma)^{-1}(E^{-1}(U))
\ \subset\  E^{-1}(\bigcap_{\gamma\in\G}\rho(\gamma)^{-1}(U))\  = 
\  {\rm Ker}(E).$$
Since the kernel of $E$ is discrete, there exists
an open set $V\subset E^{-1}(U)$  such that 
$V\cap {\rm ker}(E) = \{0\}$.
From the above identity it is easy to see that $V$ is an
expansive neighborhood of $0$
for the action $\rho_{e}$. Now the given assertion follows from
Proposition \ref{p1}.
$\hfill \Box$
\subsection{Metrics induced by endomorphism actions}
Let  $\delta$ be the function from
${{\mathbb{S}}^{1}}$ to $\R$ defined by
$\delta(z) = \mbox{inf}\{ |t|\ |\ t\in{\R},\ e^{2\pi i t} = z\}.$
For any solenoid $G$ and  $A\subset \widehat{G}$ 
we define functions \, $d_{A}: G \times G\rightarrow [0,1]$
and $d_{A}: L(G) \times L(G)\rightarrow [0,\infty ]$ by

$$\begin{array}{lll}
\vspace*{2mm}
 d_{A}(g,h) & = & {\sup_{\chi\in A}}\ \delta\circ\chi(g^{-1}h),\\
  d^{*}_{A}(p,q) & = & \sup_{\chi\in A}| p(\chi) - q(\chi)|.
\end{array}$$

\medskip
It is easy to see that if $A$ generates $\widehat{G}$ as an
abelian group then  $d_{A}$ is a metric on $G$ and  $d^{*}_{A}$
is a metric on the subspace $L(G)_{A}\subset L(G)$ consisting
of all $p$ in $L(G)$ with $d^{*}_{A}(0,p) < \infty$. 
For $A\subset \widehat{G}$ and $r > 0 $, we define
 open sets $B_{A}(r)\subset G$ and $B^{*}_{A}(r)\subset L(G)$ by

$$ \begin{array}{lll}
\vspace*{2mm}
B_{A}(r) & = & \{ g\in G\ |\  d_{A}(e,g)< r \},\\ 
B^{*}_{A}(r) & = & \{ p\in L(G)\ |\ d^{*}_{A}(0, p) < r \}. 
\end{array}$$

\medskip
Suppose $\rho$ is an endomorphism action of a semigroup $\Gamma$ on
a solenoid $G$ and  $A$ is a $\rho$-basis of $\widehat{G}$. Our
next proposition shows that expansiveness of $\rho$ and existence of
 non-trivial bounded $\rho$-orbits in $L(G)$ can be characterized in terms 
of the metrics 
$d_{A}$ and $d^{*}_{A}$ respectively.
\begin{proposition}\label{bdv}
Let $\G$ be a semigroup and $\rho$ be an endomorphism action of
$\G$ on a solenoid $G$. Let $A$ be a $\rho$-basis of $\widehat{G}$.
 Then we have the following.
\begin{enumerate}
\item{ $\rho$ is expansive  if and only if
$B_{A}(\epsilon) = \{ e\}$ for some $\epsilon > 0$.}
\item{  $\rho_{e}$ has a nontrivial bounded
orbit if and only if 
$B^{*}_{A}( C) \neq \{ 0\}$ for some $C > 0$.}
\end{enumerate}
\end{proposition}
\proof
1) Suppose $B_{A}(\epsilon) = \{ e\}$ for some $\epsilon > 0$. We
choose a finite set $F\subset A$ such that 
$A = \bigcup_{\gamma}{\widehat{\rho}}(\gamma)(F)$. Note that
for any $g$ in $G$, the $\rho$-orbit of $g$ is contained in
$B_{F}(\epsilon)$  only if $g\in B_{A}(\epsilon)$. Hence
 $B_{F}(\epsilon)\subset G$ is an expansive 
neighborhood of $e$ for the action $\rho$.

On the other hand suppose $U\subset G$ is  an expansive 
neighborhood of $e$. We
choose a finite set $F\subset A$ such that 
$A = \bigcup_{\gamma}{\widehat{\rho}}(\gamma)(F)$. 
Since $A$ generates
$\widehat{G}$ as an abelian group, it follows that 
 $B_{F}(\epsilon)\subset U$  for some $\epsilon > 0$.
Now it is easy to see that
$$ B_{A}(\epsilon)\, =
\, \bigcap_{\gamma}\rho(\gamma)^{-1}(B_{F}(\epsilon))\,
\, \subset \bigcap_{\gamma}\rho(\gamma)^{-1}(U)\, = \, \{ e\}. $$
2)  Suppose $p$ is a non-zero element of $L(G)$ such that
the $\rho_{e}$-orbit of $p$ is  bounded.
 We fix an element $\phi$ in
$ {\widehat{G}}$ and observe that
$$ \{ q(\phi)\ |\ q \in\rho_{e}(\Gamma)(p)\}
=  \{ p(\psi)\ |\ \psi\in{\widehat{\rho}}(\Gamma)(\phi)\}.$$
Since $q\mapsto q(\phi)$ is a continuous map from $L(G)$ to ${\mathbb R}$
and the $\rho_{e}$-orbit of $p$ is  bounded, 
the right hand side is a bounded subset of $\R$.
Hence  $p$, restricted to the $\widehat{\rho}$-orbit of
$\phi$ is bounded. Since $\phi$ is arbitrary and
 $A$ is a union of finitely many
$\widehat{\rho}$-orbits, we conclude that $p|_{A}$ is bounded i.e. 
$p\in  B^{*}_{A}(C)$ for some $C > 0$.

Conversely, suppose $p$ is an element of
$ B^{*}_{A}(C)$ for some $C > 0$. We fix an element $\phi$ in
$\widehat{G}$. Since $A$
generates $\widehat{G}$ as an abelian group, there exists a positive 
integer $l$ and
$a_{1},\ldots ,a_{l}$ in $A$ 
such that 
$\phi = a_{1} + \cdots + a_{l}$. Since $A$ is invariant under
$\widehat{\rho}$,
every element in the $\widehat{\rho}$-orbit of $\phi$ can be
written as sum of $l$ elements of $A$.
Hence
 $|p(\psi)| < lC$ for all $\psi$ in the $\widehat{\rho}$-orbit of
$\phi$ i.e. $|q(\phi)| < lC$ for all $q$ in the $\rho_{e}$-orbit
of $p$. Since $\phi$ is arbitrary, this implies that the 
$\rho_{e}$-orbit
of $p$ is bounded.
$\hfill \Box$
\subsection{Actions of finitely generated semigroups}
As remarked earlier, for an arbitrary solenoid $G$,
 the map $E : L(G)\rightarrow G$
need not be surjective.
However we will prove the following result which will be a crucial step
in the proof of Theorem B.
\begin{theorem}\label{main}
Let $\G$ be a finitely generated semigroup and $\rho$ be an endomorphism
action of $\G$
on a solenoid $G$. Let $A$ be a $\rho$-basis of $\widehat{G}$.
Then for any $C > 0$ there exists an $\epsilon > 0$ such that
$B_{A}(\epsilon) \subset E(B^{*}_{A}(C))$. 
\end{theorem}
Note that this theorem, together with Propositions  
\ref{basis}, \ref{div} and \ref{bdv} implies Theorem B when $\Gamma$ is 
finitely generated.
\rem
Earlier we have observed that  $E$ is a continuous homomorphism 
and the kernel of $E$ is
discrete. Hence the above theorem implies that the map $E$,
 restricted to $ L(G)_{A}$,
is a local homeomorphism if $L(G)_{A}$ and $G$ are equipped with the
topologies induced by  $d^{*}_{A}$ and $d_{A}$ respectively.

\bigskip
Before beginning the proof of Theorem \ref{main} we introduce a few
notations.

\medskip
Let $H$ be a discrete abelian group and $A$ is a finite subset of
$H$. Then  by ${\bf Q}(A)$ we denote the set of all
$h\in{H}$ which satisfies a linear equation of the form 
$$ n_{0}h = \sum_{j=1}^{r} n_{j}a_{j},$$
where $a_{1},\ldots ,a_{r}$ are elements of $A$
and $n_{0},n_{1},\ldots ,n_{r}$
are integers with $n_{0} \ne 0$.
Also for all $k > 0$, by ${\bf Q}_{k}(A)$  we denote the set of all
$h\in{H}$ which satisfies a linear equation as above
with
$|n_{0}| + |n_{1}| + \cdots + |n_{r}| \le k$.
A set $A \subset H$ is said to be a {\it $k$-regular \/}
if  there exists an increasing sequence of finite sets
$A_{1}\subset A_{2}\subset\cdots \subset A$ satisfying
$$ \bigcup_{i = 1}^{\infty} A_{i} = A,\ \ A_{n}\subset {\bf Q}_{k}(A_{n-1})
\ \ \forall n\ge 2.$$

\smallskip
Now we show that if $\rho$ is an  endomorphism action of a 
finitely generated semigroup on a solenoid $G$, then any 
$\rho$-basis is $k$-regular for some $k$.
\begin{proposition}\label{GF}
Let $\G$ be a finitely generated semigroup and $\rho$ be an
endomorphism action of $\G$
on a solenoid $G$. Then any $\rho$-basis of $\widehat{G}$ is 
$k$-regular for some $k > 0$.
\end{proposition}
\proof
Let $A$ be a $\rho$-basis of $\widehat{G}$. Since $A$ generates 
$\widehat{G}$ as an abelian group and $\widehat{G}$ has finite rank,
 there exists a finite set 
$F_{1}\subset A$ such that ${\bf Q}(F_{1}) = \widehat{G}$.
Also,  $A = \bigcup_{\gamma}{\widehat{\rho}}(\gamma)(F_{2})$ for some
finite set $F_{2}\subset A$. We define
$F = F_{1}\cup F_{2}$. 
 Let $S\subset \Gamma$
be finite set which generates $\Gamma$ as a semigroup. 
Since  ${\bf Q}(F) = \widehat{G}$, for any $\gamma$ in $S$ and $a_{0}$
 in $F$,
  ${\widehat{\rho}}(\gamma)(a_{0})$
satisfies a linear equation of the form
 $$ n_{0}\cdot{\widehat{\rho}}(\gamma)(a_{0}) = \sum_{j=1}^{r} n_{j}a_{j},$$
where $a_{1},\ldots ,a_{r}$ are elements of $F$ 
and $n_{0},n_{1},\ldots ,n_{r}$
are integers with $n_{0}\ne 0$. Since  $S$ is finite, it follows that
there exists 
a positive integer $k$ such that 
 ${\widehat{\rho}}(\gamma)(F) \subset {\bf Q}_{k}(F)$ for all $\gamma$ in $S$.
For $n\ge 1$ we define $S_{n}\subset \G$ and
$A_{n}\subset A$ by

$$ S_{n} = \{ \gamma_{1}\cdots \gamma_{n}\ |\ 
\gamma_{1},\ldots ,\gamma_{n}\in S\},\ \ 
A_{n}  =  \bigcup_{\gamma\in S_{n}}{\widehat{\rho}}(\gamma)(F).$$

Let $\gamma_{n}$ be any element of $S_{n}$. We choose 
$\gamma_{n-1}\in S_{n-1}$ and $\gamma_{0}\in S$ such that 
$\gamma_{n} = \gamma_{n-1} \cdot \gamma_{0}$.
Since ${\widehat{\rho}}(\gamma_{0})(F)\subset {\bf Q}_{k}(F)$, we see   that

$$ {\widehat{\rho}}(\gamma_{n})(F)\ \subset 
\ {\widehat{\rho}}(\gamma_{n-1})({\bf Q}_{k}(F))
\ \subset\  {\bf Q}_{k}({\widehat{\rho}}(\gamma_{n-1})(F)).$$

\medskip
\noindent
This shows that $A_{n}\subset {\bf Q}_{k}(A_{n-1})$ for all $n\ge 2$.
Since $S$ generates $\Gamma$, it is easy to see that 
$\bigcup A_{i} = {\widehat{\rho}}(\Gamma)(F) = A$.
$\hfill \Box$

\bigskip
Now we turn to the proof of Theorem \ref{main}.
We will use the  following lemma.
\begin{lemma}\label{fin}
Let $G$ be a solenoid
and $F$ be a finite subset of ${\widehat{G}}$. Then for any $C > 0$ there exists
an $\epsilon > 0$ such that for any $g\in B_{F}(\epsilon) $ there exists
$p$ in $ B^{*}_{F}(C)$ satisfying $\phi(g) = \phi\circ E(p)$ for all
$\phi$ in $F$.
\end{lemma}
\proof
Let $\widehat{H}\subset \widehat{G}$ be the subgroup
generated by $F$ and let $H$ be the dual of $\widehat{H}$.
Since $\R$ is divisible, any homomorphism from $\widehat{H}$ to $\R$
can be extended to a homomorphism from ${\widehat{G}}$ to $\R$. Therefore
it is enough to consider the case when $G = H$.
In that case there exists $n\ge 1$ and an isomorphism
 $\theta : G\rightarrow{\mathbb T}^{n}$
 such that ${\rm exp} \circ {\rm d}\theta = \theta \circ E $.
Since $B^{*}_{F}(C)$ is an open subset of $L(G)$ and
 ${\rm exp }: L({\mathbb T}^{n})\rightarrow {\mathbb T}^{n}$ is a local 
homeomorphism, it is easy to see that
$E(B^{*}_{F}(C))$ contains an open neighborhood of $e$ in $G$.
Since $F$ generates $\widehat{G}$,  there exists an
$\epsilon > 0$ such that
$ B_{F}(\epsilon) \subset E(B^{*}_{F}(C))$.
 This proves the lemma.
$\hfill \Box$

\bigskip
{\bf Proof of Theorem \ref{main} }
Suppose $A$ is a $\rho$-basis of $\widehat{G}$. By Proposition  
\ref{GF} $A$ is $k$-regular for some $k > 0$. 
Clearly  without loss of generality we may assume that
$C < 1/k$.
Let $A_{1}\subset A_{2}\subset\cdots \subset A$ be an increasing
sequence of finite sets such that
$$ \bigcup_{i } A_{i} = A,\ \ A_{n}\subset {\bf Q}_{k}(A_{n-1})
  \ \ \forall n\ge 2. $$ 
 Applying  Lemma \ref{fin} we choose a positive $\epsilon < C$ 
such that for all $h\in B_{A_{1}}(\epsilon)$ there exists a 
$q\in B^{*}_{A_{1}}(C)$ satisfying $\chi(E(q)) = \chi(h)$ for all
$\chi$ in $A_{1}$.
Let $g$ be an element of $B_{A}(\epsilon)$. We choose
a $p$ in $B^{*}_{A_{1}}(C)$ such that $\chi(E(p)) = \chi(g)$ for all
$\chi$ in $A_{1}$. We claim that

$$|p(\chi)| < C,\  \chi(E(p)) = \chi(g)\ \ \ \forall 
\chi\in A.$$

 By our choice of $p$ this is true for all $\chi$ in
$A_{1}$.
Suppose this holds for all $\chi$ in $A_{n}$. 
Let $a_{0}$ be an element of $A_{n+1}$.
We choose integers $n_{0},n_{1},\ldots ,n_{r}$ and 
$a_{1},\ldots ,a_{r}\in A_{n}$ such that $n_{0} \ne 0$,
$|n_{0}| + |n_{1}| + \cdots + |n_{r}| \le k$ and
$$ n_{0}a_{0} = \sum_{j=1}^{r} n_{j}a_{j}.$$
Since  $g\in B_{A}(\epsilon)$,
 there exists $\alpha\in \R$ such that
$|\alpha | < \epsilon < C$ and $e^{2\pi i \alpha} =  a_{0}(g)$.
Now from the above equation it is easy to see that
$$e^{2\pi i n_{0}\alpha}\  = \  a_{0}(g)^{n_{0}}\ 
= \ \prod_{1}^{r} a_{i}(g)^{n_{i}}.$$
Since $a_{1},\ldots ,a_{r}\in A_{n}$, it follows that
$$ e^{2\pi i n_{0}p(a_{0})} =
 \prod_{1}^{r}e^{2\pi i n_{j}p(a_{j})} = 
\prod_{1}^{r}a_{i}\circ E(p)^{n_{i}}
= \prod_{1}^{r} a_{i}(g)^{n_{i}}.$$
Therefore $e^{2\pi i n_{0}(\alpha - p(a_{0}))} = 0$ i.e.
$n_{0}(\alpha - p(a_{0}))\in \Z$. On the other hand

$$|n_{0}(\alpha - p(a_{0}))|\, \le\,  
 n_{0}|\alpha| + |p(n_{0}a_{0})| \, \le  n_{0}C + \sum_{1}^{r}
n_{i}|p(a_{i})|.$$

\smallskip
\noindent
Now applying our induction hypothesis we see that 

$$|n_{0}(\alpha - p(a_{0}))|\, \le\,  
 C\sum_{j = 1}^{r} |n_{j}| \, < \,  Ck\, < \, 1.$$
\medskip
\noindent
Hence $\alpha - p(a_{0}) = 0$, which implies that
$|p(a_{0})| = |\alpha | < C$ and 
$a_{0}(E(p)) = e^{2\pi i \alpha} =  a_{0}(g)$. This proves the claim.
Since $A$ generates $\widehat{G}$ as an abelian group, from the above
claim
we deduce that for any $g\in B_{A}(\epsilon)$ there exists a $p\in
B^{*}_{A}(C)$ such that $E(p) = g$.
This  completes the proof.
$\hfill \Box$
\subsection{Proof of Theorem B}
Now we prove Theorem B for an endomorphism action $\rho$ of
an arbitrary semigroup $\Gamma$ on a solenoid $G$. 
If $(G,\rho)$ is expansive then from Proposition
 \ref{basis} and Proposition \ref{div}
 it follows that a) and b) are satisfied. Conversely, suppose
the conditions a) and b) are satisfied. 
By Lemma \ref{key}
there exists a finitely generated subsemigroup $\Gamma_{0}$ of $\G$
 such that for all $p$ in $L(G)$, the $\Gamma_{0}$-orbit of
$p$ under the action $\rho_{e}$ is unbounded. We construct an
endomorphism action
$(H,\sigma)$ of $\Gamma_{0}$ as follows :

Since $\widehat{G}$ has
finite rank,  
there exists a finite set  $F_{1}\subset \widehat{G}$ 
 such that
${\bf Q}(F_{1}) = \widehat{G}$. Also, from a) it follows that there exists
a finite set  $F_{2}\subset \widehat{G}$ 
 such that the set  
$ \bigcup_{\gamma}{\widehat{\rho}}(\gamma)(F_{2})$ 
generates $\widehat{G}$ as
an abelian group. We define $F = F_{1}\cup F_{2}$.
Let ${\widehat{H}}\subset {\widehat{G}}$ be the subgroup generated
by the set 
$$ A_{0} = \bigcup_{\gamma\in\Gamma_{0}}{\widehat{\rho}}(\gamma)(F)$$
and let $\widehat{\sigma}$ be the $\Gamma_{0}$-action on
${\widehat{H}}$
defined by ${\widehat{\sigma}}(\gamma) = 
{\widehat{\rho}}(\gamma)|_{\widehat{H}}$.
Let $H$ denote the dual of ${\widehat{H}}$ and let $\sigma$ denote the
endomorphism
action of $\Gamma_{0}$ on $H$ which is the dual of
 $\widehat{\sigma}$.

\medskip
{\bf Claim :} $(H,\sigma)$ is expansive.

\medskip
Since  $\widehat{G} =  {\bf Q}(F)\subset 
{\bf Q}(\widehat{H})$, it follows that the restriction map 
$p\mapsto p|_{\widehat{H}}$ is a $\Gamma_{0}$-equivariant linear
isomorphism from $L(G)$ onto $L(H)$. This implies that 
for every non-zero $q\in L(H)$, the ${\sigma}_{e}$-orbit of $q$ is
unbounded in $L(H)$. 
Also it is easy to see that $A_{0}$ is a $\sigma$-basis of 
$\widehat{H}$. Now applying Theorem \ref{main} and Proposition \ref{bdv}
we see that $(H,\sigma)$ is expansive, which proves the claim.  

\medskip
 Let $V$ be an expansive neighborhood of
$e$ in $H$. Let $i : \widehat{H}\rightarrow \widehat{G}$ be the  
inclusion map and let $\pi : G \rightarrow H$ be the dual of $i$.
We claim that $\pi^{-1}(V)$ is an expansive neighborhood for
 the action $\rho$.
To see this we choose any $g$ in $G$ such that
the $\rho$-orbit of $g$ is contained in $\pi^{-1}(V)$.  
It is easy to see that for any  $\gamma\in\Gamma_{0}$,
 $$\sigma(\gamma)\circ \pi(g) = \pi\circ\rho(\gamma)(g).$$
Therefore
the $\sigma$-orbit of $\pi(g)$ is contained in
$V$. Since $V$ is an expansive neighborhood for $\sigma$, this implies
that  $\pi(g) = e$. In particular
$\phi(g) = 1 $ for all $\phi$ in $F$. For any $\gamma$ in
$\G$, replacing $g$ by $\rho(\gamma)(g)$ and
applying the same argument we see that
$$ {\widehat{\rho}}(\gamma)(\phi)(g) = 
\phi\circ\rho(\gamma)(g) = 1 \ \ \forall\gamma\in{\G},\ \phi\in{F}.$$
Since $A = \bigcup_{\gamma \in \G}{\widehat{\rho}}(\gamma)(F)$ 
generates ${\widehat{G}}$, 
it follows that  $g =e$. This completes the proof.
$\hfill \Box$
 \rem
Let $\Gamma$ be the group of positive rational numbers under
multiplication, ${\mathbb Q}$ be the  group of
rational numbers under addition and $\widehat{\rho}$ be the
 action of $\Gamma$ on  ${\mathbb Q}$ defined by 
$\widehat{\rho}(\gamma)(x) = \gamma\cdot x$. Let $\rho$ be the dual 
action of $\Gamma$ on $\widehat{\mathbb Q}$. Then from Theorem B it
follows that $( {\widehat{\mathbb Q}}, \rho )$ is expansive.
On the other hand it is easy to see that for any finitely
generated subgroup $\Gamma_{0}\subset \Gamma$,  ${\mathbb Q}$ is not
finitely
generated as a ${\mathbb Z}(\Gamma_{0})$-module. Hence no finitely
generated subgroup $\Gamma_{0}\subset \Gamma$ acts expansively
on $ {\widehat{\mathbb Q}}$ under the action $\rho$. Since
$\Gamma$ is abelian, this shows that analogues of Corollary \ref{im}
and Corollary \ref{nil} do not hold when $G$ is  a solenoid.

\bigskip
\noindent
Address : School of Mathematics, Tata Institute of Fundamental Research,
Mumbai 400005, India.\\
E-mail : siddhart@math.tifr.res.in

\end{document}